\documentclass[10pt,a4paper,twoside]{article}
\usepackage{amsfonts,amssymb}

\font\Bigtit=cmr10 scaled \magstep 4
\font\ebf=cmbx8
\font\erm=cmr8
\font\eit=cmti8

\usepackage{graphicx}

\begin{document}

\thispagestyle{empty}

\begin{flushright}
PL ISSN 0459-6854
\end{flushright}
\vspace{0.5cm}
\centerline{\Bigtit B U L L E T I N}
\vspace{0.5cm}
\centerline{DE \ \  LA \ \  SOCI\'ET\'E \ \  DES \ \  SCIENCES \ \ ET \ \ DES \
\ \ LETTRES \ \ DE \ \ \L \'OD\'Z}
\vspace{0.3cm}
\noindent 2002\hfill Vol. LII
\vspace{0.3cm}
\hrule
\vspace{5pt}
\noindent Recherches sur les d\'eformations \hfill Vol. XXXVI
\vspace{5pt}
\hrule
\vspace{0.3cm}
\noindent pp.~45--65

\vspace{0.6cm}

\noindent {\it Andrzej~K.~Kwa\'sniewski}

\vspace{0.5cm}

\noindent {\bf ON SIMPLE CHARACTERISATIONS OF SHEFFER\\ $\psi$-POLYNOMIALS AND
RELATED PROPOSITIONS\\ OF THE CALCULUS OF SEQUENCES}

\vspace{0.5cm}

\noindent {\ebf Summary}

{\small  A ``Calculus of Sequences" had started by the 1936 publication of
Ward suggesting the possible range for extensions of  operator calculus of
Rota-Mullin, considered by several authors and after Ward. Because of
convenience we shall call the Wards calculus of sequences in its afterwards
elaborated form -- a $\psi$-calculus. The notation used by Ward, Viskov,
Markowsky and Roman is accommodated in conformity with Rota's way of
exposition. In this manner $\psi$-calculus becomes in parts almost
automatic extension of finite operator calculus. The $\psi$-extension
relies upon the notion of $\partial _{\psi}$-shift invariance of operators.
At the same time this calculus is an example of the algebraization of the
analysis -- here restricted to the algebra formal series. The efficiency of
the notation used is further exemplified  among others by easy proving of
some Sheffer $\psi$-polynomials characterisation theorems as well as
Spectral Theorem. $\psi$-calculus results may be extended to Markowsky
``$Q$-umbral calculus", where $Q$  stands for a generalised difference
operator (not necessarily $\partial _{\psi}$-shift invariant)  i.e. the one
lowering the degree of any polynomial by one.}

\vspace{0.4cm}

\noindent {\ebf Contents}

{\small 1.~Introduction. -- 2.~Primary definitions, notation  and
the general observations. -- 3.~The general picture of  $End(P)$. --
4.~Characterisations of Sheffer $\psi$-polynomials and related
propositions. -- 5.~Miscellaneous remarks and indications of several
applications.}


\renewcommand{\thesubsection}{\arabic{subsection}.}
\subsection{Introduction}
\setcounter{equation}{0}

Already  since seventies of the past century it had been realised that the
notions and methods of Rota-Mullin finite operator calculus (started in [1]
and developed in [2, 3] -- see also  [4--6]) might be extended to the use
of {\it any} polynomial sequences  $\{ p_{n}\}_{0}^{\infty}\ ({\rm deg}\ p_{n}=n)
$ instead of those of binomial type only. The foundations of such an
extension were led in 1975  by Viskov in  [7] and then developed further
in [8] \linebreak -- with very substantial reference to Boas and  Buck [9, 10]  --
however without reference to Ward paper [11].

The decisive contribution to this conviction  (on the possibility to
develop umbral calculus for  {\it any} polynomial sequences
$\{p_{n}\}_{0}^{\infty}$ instead of those of binomial type only) we owe to
Markowsky [12] who apparently was not acquainted with Viskov papers [7, 8]
at that time and neither [7]  nor Ward paper [11] are quoted in [12].
Markowsky paper dealt with -- as we would call it now [13] -- generalised
Sheffer polynomials.

One may learn much more from later  works on subsequent and recent
progress in extending the operator methods onto  ``calculus of sequences".
The reader is referred to the up-date position [14] by Loeb. This  is
extensive and exhaustive source of references and the survey of widespread
developments  over past decades. With such enriched experience -- any  way
--  we are led back to the Ward conception  of calculus of sequences  as
the source.

In [15, 16] special type generalised Sheffer polynomials are called
Sheffer \linebreak $\psi$-polynomials in order to accommodate the author's
Rota-oriented notation with  the notation used by Viskov in [7]  on one
side and notation of Rota  in [3] on the other side. This (hopefully well
aimed) notation enables to adopt from [3, 12]  formulations and methods for
proving majority of  ``$\psi$-propositions" or  ``$Q$-propositions" (see
definitions 2.5. and 2.6.) in the general  case of  operator difference
calculus taking care of  and related to   {\it any} polynomial sequences
$\{p_{n}\}_{0}^{\infty}$ instead of those of binomial type only  or
$\psi$-binomial type only.  As an example of such presentation
characterisation of Sheffer $\psi$-polynomials theorems and Spectral
Theorem [16] are given among others.  The note is organised as follows.
Firstly the presentation of general picture of algebra of linear operators
on polynomial algebra is accomplished. Then  characterisations of  Sheffer
$\psi$-polynomials and related propositions of the calculus of sequences follow.
At the end link is given to specific formulation of  $q$-umbral calculus by
Cigler [17] and Kirschenhofer [18]. This formulation  might be related -- as
noted in  [15] -- to the so-called quantum groups [19]. The relevant
$q$-extensions of  what is now called sometimes -- Rota-Mullin calculus --
we owe also to earlier authors such as Al.-Salam, Carlitz  (see  Chihara
[20, 21] ), Goldmann, Rota,  Andrews,  Ihrig, Ismail, Verna and others.
For the corresponding references see [22] and Roman's book [23].  There --
in chapter 6  one finds also the first principal results on
{\it ``nonclassical umbral calculi"}. (It seems so that the author of
[19] was  not acquainted with Viskov and Markowsky papers at that time). In
[23] Roman refers to Ward's papers including [11], which was not the case
in Roman's early publications on the subject. Meanwhile in  [11] Ward
proposes a generalisation of a large portion of calculus of finite
differences almost in the spirit of later  finite operator calculus. Ward
accurately and relevantly called this generalisation -- the calculus of
sequences. (As a matter of fact -- the $q$-calculus of sequences is already
present as an example in Ward`s publication [11]). Here we shall deal with
the past century further development of the finite operator calculus
formulation of this calculus of sequences. To this end  we choose  as a
motto for activities of the type to be presented -- a quotation from
Roman's book (Chapter~6, p.~162  in [23] -- with notation and reference
changed into the one used in this note): {\it ``Let $n_{\psi}! \equiv
c_{n}$  be a sequence of nonzero constants. If $n!$ is replaced by $c_{n}$
throughout the preceding theory, then virtually all of the results remain
true, mutatis mutandis. In this way each sequence  $c_{n}$ gives rise  to a
distinct umbral calculus}.

{\it Actually, Ward {\rm [11]}  seems to have been the first to suggest
such a generalisation  (of the calculus of finite differences)  in  1936,
but the idea remained relatively undeveloped until quite recently, perhaps
due to a feeling that it was mainly generalisation for its own  sake. Our
purpose here is to indicate that this is not the case."}


\subsection{Primary definitions, notation  and the general
observations}
\setcounter{equation}{0}

In the following we shall consider the algebra $P$ of polynomials  $P =
{\rm \bf F}[x]$ over the field  {\bf F}  of characteristic zero. All
operators or functionals studied here are to be understood as {\it linear}
operators or functionals  on  $P$. It shall be easy to see that they are
always well defined.  Whenever we say polynomial sequence $\{p_{n} \}
_{0}^{\infty}$ we shall always mean such sequence  $\{p_{n} \}
_{0}^{\infty}$  of polynomials that
$\deg\  p_{n} = n$.

Let $s = \{ s_{n} \}_{n \geq 1}$   be  an arbitrary sequence of numbers
such that $s_{n} \neq 0$, $n \in N$. $s$-binomial coefficients  are then
defined with help of the generalised factorial
$n_{s}!=s_{1}\ s_{2}\ s_{3}\ ...\ s_{n}$ and
$n_{s}^{\underline{k}}=n_{s}(n-1)_{s}...(n-k+1)_{s}$ in a usual way as
$\left(\begin{array}{c} n\\k  \end{array} \right) \equiv
\frac{n_{s}^{\underline{k}}}{k_{s}!}$.
Consider  $\Im$ -- the family of functions' sequences (in conformity with
Viskov notation) such that
\begin{eqnarray}
\Im &=& \{\psi ;\ R \supset [a,b];\ q \in[a,b];\ \psi (q):Z \to F;\nonumber
\\
&& \psi _{0}(q)=1;\ \psi _{n}(q) \neq 0;\ \psi _{-n}(q)=0;\ n \in N\}
\nonumber \end{eqnarray}
\noindent Following Roman [23--25] we shall
call  $\psi=\{\psi _{n}(q) \}_{n\geq 0}$; $\psi _{n}(q) \neq 0$; $n \geq 0$
and $\psi _{0}(q)=1$ an {\it admissible sequence}. Let now $n_{\psi}$ --
(in conformity with Viskov notation) -- denotes  [16, 15]
\[
n_{\psi} \equiv \psi _{n-1}(q) \psi_{n}^{-1}(q).
\]

\noindent Then we have for the $\psi$-factorial $n_{\psi}! \equiv \psi
_{n}^{-1}(q)\equiv
n_{\psi}(n-1)_{\psi}(n-2)_{\psi}(n-3)_{\psi}...2_{\psi}1_{\psi}$;
$0_{\psi}!=1$ and $\left(\begin{array}{c} n\\k \end{array}  \right)_{\psi}
\equiv \frac{n_{\psi}^{\underline{k}}}{k_{\psi}!}$
for   $\psi$-binomial coefficients where
$n_{\psi}^{\underline{k}}=n_{\psi}(n-1)_{\psi}...(n-k+1)_{\psi}$ while
$\exp _{\psi}\{y \}=\sum _{k=0}^{\infty} \frac{y^{k}}{k_{\psi}!}$ defines
the $\psi$-exponential series.
In this note we need not to specify any additional conditions on
$\psi_{n}(q)$'s which would lead us to specifications of the
$\psi$-calculus as for example the $q$-umbral calculus which is obtained
with the following choice of an admissible $\psi$:
\[
\psi _{n}(q)= \frac{1}{R(q^{n})!},\quad R(x)=\frac{1-x}{1-q}.
\]

\vspace{3mm}

\noindent {\it Definition 2.1.}  Let  $\psi$  be admissible. Let
$\partial _{\psi}$ be the linear  operator lowering degree of polynomials by
one defined according to $\partial _{\psi}x^{n}=n_{\psi}x^{n-1}$;  $n \geq
0$. Then $\partial _{\psi}$ is called the {\it$\psi$-derivative.}

\vspace{3mm}

\noindent {\it Remark 2.1.} The choice $\psi _{n}(q)=[R(q^{n})!]^{-1}$  and
$R(x)=\frac{1-x}{1-q}$ results in the well known $q$-factorial
$n_{q}!=n_{q}(n-1)_{q}!;$ $1_{q}!=0_{q}!=1$ while the $\psi$-derivative
$\partial _{\psi}$ becomes now the Jackson's derivative [16, 15]
$\partial _{q}$:
\[
\left( \partial _{q} \varphi \right)\left( x \right)= \frac {\varphi
(x)-\varphi(qx)}{(1-q)x}. \]

\vspace{3mm}

\noindent {\it Example 2.1.} [11]   Let $\{n_{\psi}\}_{n \geq 0}$ be the
solution of a linear recurrence of  {\it r}-th order for which its
characteristic polynomials has  $r$   distinct roots  $\{ \alpha _{k}
\}_{1}^{r}$. Then the general solution reads:  $n_{\psi}=\beta _{1}\alpha
_{1}^{n}+\beta _{2}\alpha _{2}^{n}+...+\beta _{r}\alpha _{r}^{n}$ and one
now should impose the two conditions: $\sum_{k=1}^{r}\beta _{k}=0$ (because
$0_{\psi}=0$) and $\sum_{k=1}^{r}\beta _{k}\alpha _{k}=1$ (because
$1_{\psi}=1$). Naturally for any formal series $F(x)=\sum _{k \geq
0}f_{k}x^{k}$, $\partial _{\psi}F(x)=\sum_{k \geq 0}k_{\psi}f_{k}x^{k-1}$
therefore $\partial _{\psi}F(x)=\frac{1}{x} \sum _{k=1}^{r}\beta
_{k}F\left( \alpha _{k}x \right)$. For $r = 2$ and $\alpha _{1}=q$,
$\alpha _{2}=1$, $\beta _{1}=(q-1)^{-1}$, $\beta _{2}=-\beta _{1}$ the
$\psi$-derivative $\partial _{\psi}$ becomes Jackson's derivative i.e.
$\partial _{\psi}=\partial _{q}$ [11, 16]. Another specific choice is the
Fibonacci sequence
\[
F_{n}=\frac {1}{\sqrt{5}} \left( \frac{1+\sqrt{5}}{2}
\right)-\frac{1}{\sqrt{5}}\left( \frac{1-\sqrt{5}}{2}
\right);\qquad n\geq 0.  \]
\noindent Then $\psi$-derivative becomes ``Fibonacci"
derivative i.e. $\partial _{\psi}=\partial _{F}$ with $n_{\psi} \equiv
n_{F}\equiv F_{n}$. The quite exceptional property of $\left(
\begin{array}{c} n\\k \end{array} \right)_{\psi} \equiv
\frac{n_{\psi}^{\underline{k}}}{k_{\psi}!}$ $\psi$-binomial coefficients is
that $\left(
\begin{array}{c} n\\k \end{array} \right)_{F} \equiv
\frac{n_{F}^{\underline{k}}}{k_{F}!}$ are integer numbers (see chapter~6,
exercises~86 in [26]) in the case of ``Fibonomial" coefficients.

\vspace{3mm}

\noindent {\it Definition 2.2.} Let  $E^{y}(\partial _{\psi})\equiv
\exp _{\psi} \{ y \partial _{\psi}\}=\sum _{k=0}^{\infty}\frac
{y^{k}\partial _{\psi}^{k}}{k_{\psi}!}$. $E^{y}(\partial _{\psi}) $ is
called the {\it generalised translation operator.}

\vspace{3mm}

\noindent {\it Note 2.1.} [16, 15]  $E^{a}(\partial _{\psi})\ f(x)\equiv
f(x+ _{\psi}a)$; $(x+ _{\psi}a)^{n} \equiv E^{a} (\partial _{\psi})x^{n}$;
$E^{a}(\partial _{\psi})f=\sum _{n \neq 0}\frac{a^{n}}{n_{\psi}!} \partial
_{\psi}^{n}f$;  and in general $(x + _{\psi} a)^{n} \neq (x +
_{\psi}a)^{n-1}(x + _{\psi} a)$. Note also that in general $(1 + _{\psi}
(-1))^{2n+1} \neq 0$; $n \geq 0$ though $(1 + _{\psi} (-1))^{2n} = 0$; $n
\geq 1$. We learn from {\rm [11]} the following.

\vspace{3mm}

\noindent {\it Note 2.2.}  $\exp _{\psi} (x + _{\psi}y) \equiv \exp
_{\psi} \{x\} \exp _{\psi} \{y\}$  while in general $\exp _{\psi}\{
x+y\} \neq \exp _{\psi}\{ x \}\exp_{\psi} \{y\}$. Possible consequent
utilisation of the identity $\exp _{\psi}(x+ _{\psi}y) \equiv \exp _{\psi}
\{x\}\exp _{\psi} \{y\}$ is quite encouraging. It leads among others to
``$\psi$-trigonometry" either $\psi$-elliptic or $\psi$-hyperbolic via
introducing $\cos _{\psi}$, $\sin _{\psi}$ {\rm [11]}, $\cosh _{\psi}$,
$\sinh _{\psi}$ or in general $\psi$-hyperbolic functions of m-th order
$\left\{h_{j}^{(\psi)} (\alpha) \right\} _{j \in {\rm \bf Z}_{m}}$ defined
according to [27, 28]
\[ R \backepsilon \alpha \to h_{j} (\alpha)= \frac
{1}{m} \sum _{k \in Z_{m}}\omega ^{-kj} \exp _{\psi} \left( \omega ^{k}
\alpha \right); \quad j \in {\rm \bf Z}_{m},\ \omega =\exp \left(i \frac {2
\pi}{m}  \right). \]

\noindent where $1 <m \in N$ and ${\rm \bf Z}_{m} = \{0,1,...,m-1\}$.
 However note that elements of  $R$  (or any other field of
characteristic zero chosen) are subjected to  $\psi$-addition and
$\psi$-subtraction  {\rm [11, 16]} and that  $(x+ _{\psi}a)^{n} \neq (x+
_{\psi}a)^{n-1} (x+ _{\psi}a)$ in general.

In the following the notion of
$\partial _{\psi}$-shift invariance of operators and the notion of  a
polynomial sequence $\{p_{n} \}_{0}^{\infty}$ is of $\psi$-{\it binomial}
type are to play the crucial role.

\vspace{3mm}

\noindent {\it Definition 2.3.}  Let us denote by $End(P)$ the algebra
of all linear operators  \linebreak acting on the algebra  P of
polynomials.  Let $\sum_{\psi}=\{ T \in End(P)$; $\forall\  \alpha \in F$;
$\left[ T;\ E^{\alpha}(\partial _{\psi})\right]=0  \}$. Then $\sum _{\psi}$
is a commutative subalgebra of End (P) of F-linear operators. We shall call
these operators T: {\it $\partial _{\psi}$-shift invariant operators.}

\vspace{3mm}
\noindent {\it Definition 2.4.}   A polynomial sequence $\{p_{n}
\}_{0}^{\infty}$ is of {\it $\psi$-binomial type} if it satisfies the recurrence
\[
E^{y}\left( \partial _{\psi} \right)p_{n}(x) \equiv p_{n}\left( x+ _{\psi}y
\right)\equiv \sum _{k \geq 0} \left( \begin{array}{c} n\\k \end{array}
\right)_{\psi} p_{k}(x)p_{n-k}(y). \]

\noindent Let us express now two characterisations [7] of
polynomial sequences of $\psi$-{\it binomial} type in appealing
Rota-oriented notation as a matter of this appealingness illustration.

\vspace{3mm}

{\bf Illustration 2.1}
\begin{enumerate}
\item (Proposition 8  in  [7])  A polynomial sequence $\{ p_{n}
\}_{0}^{\infty}$ is Sheffer $\psi$-polynomial if  and only if   its
``$\psi$-generating function" is of the form:
\begin{equation}
\sum _{n \geq0} \psi _{n}p_{n}(x)z^{n}=A(z)\exp _{\psi} \left( xg(z)
\right); \end{equation}
\begin{equation}
\psi (z)=\sum_{n \geq 0}\psi _{n}z^{n};\quad \psi _{n} \neq 0; \quad
n=0,1,2,... \end{equation}

\noindent where $A(z)$, $g(z)/z$  are formal series with constant terms
different from zero.

\item Also in  [7] (Proposition 4) Viskov have proved  that polynomial
sequence $ \{p_{n} \}_{0}^{\infty} $ is  of $\psi$-binomial type if  and
only if its ``$\psi$-generating function" is of the form:
\begin{equation}
\sum_{n \geq 0} \psi _{n}p_{n} (x)z^{n}=\exp _{\psi} \left( xg (z) \right)
\end{equation}

\noindent for formal series  $g$  inverse to appropriate formal series.
Note that  for   $\psi_{n}(q)=[n_{q}!]^{-1}$  in  (2.2)
$\psi(z)=\exp_{q}\{ z \}$ and ``$\exp_{q}$ generating function" (2.3) takes
the ``$q$-umbral" form
\begin{equation}
\sum _{n \geq 0}\frac {z^{n}}{n_{q}!}p_{n}(x)=\exp _{q} \left( xg(z) \right)
\end{equation}

\noindent Polynomial sequences of $\psi$-binomial type [2, 3]  are known to
correspond in one-to-one manner to special generalised differential
operators $Q$, namely to those $Q=Q(\partial _{\psi})$  which are
$\partial _{\psi}$-shift invariant operators [7, 16, 15]. We shall deal in
this note mostly with this special case i.e. with the so-called
$\psi$-umbral calculus [16, 15].  However before to proceed let us deliver
an elementary, basic information referring to the general case of
``$Q$-umbral calculus" -- as we call it -- started by Markowsky in [12] .
\end{enumerate}

\vspace{3mm}

\noindent {\it Definition 2.5.} Let $P = F[x]$.  Let  $Q$  be a linear
map $Q: P \to P$  such that $\forall\ p \in P\ \deg (Qp)=(\deg p)-1$
$(\deg p=-1$ means $p={\rm const}= 0)$. Then $Q$ is called a {\it generalised
difference-tial operator} {\rm [12]} or {\it Gel'fond-Leontiev {\rm [7]}
operator}.  The algebra of all linear operators commuting with the
generalised difference-tial operator $Q$  is denoted by $\sum _{Q}$.

Right from the above definitions we infer  that the following
holds.

\vspace{3mm}

\noindent {\it Observation 2.1.} Let  $Q$  be as in Definition 2.5.  Let
$Qx^{n}=\sum_{k=1}^{n}b_{n,k}x^{n-k}$ where   $b_{n,1}\neq 0$. Without
loose of generality take $b_{1,1}=1$. Then $\exists \{q_{k}\}_{q\geq 2}
\subset{\rm \bf F}$ and \linebreak $\exists$ admissible $\psi$ such that
\begin{equation}
Q=\partial _{\psi}+\sum _{k \geq 2}q_{k}\partial _{\psi}^{k}
\end{equation}

\noindent if and only if    the following condition is fulfilled

\hspace{2.3cm}$b_{n,k}=\left( \begin{array}{c} n\\k \end{array}
\right)_{\psi} b_{k,k};\quad n\geq k \geq 1,\ b_{n,1}\neq 0,\
b_{1,1}=1.$\hfill(2.5$'$)

\noindent If  $\{ q_{k} \}_{q\geq 2}$ and an admissible  $\psi$  exist then
these are unique.

\vspace{3mm}

\noindent {\it Notation 2.1.} In the case (2.5$'$) is true  we shall write:
$Q=Q(\partial _{\psi})$.

\vspace{3mm}

\noindent {\it Remark 2.2.}  Note that operators of the (2.5) form
constitute a group under superposition of  formal
power series (compare with the formula (S)  in [13]).  Of course not all
generalised difference-tial operators  satisfy  (2.5) i.e. are series just
only in corresponding $\psi$-derivative   $\partial _{\psi}$ (see
Proposition 3.1). For example  [14] let
$Q=\frac{1}{2}D\hat{x}D-\frac{1}{3}D^{3}$. Then
$Qx^{n}=\frac{1}{2}n^{2}x^{n-1}-\frac{1}{3}n^{\underline{3}}x^{n-3}$.
Therefore according to Observation 2.1 $n_{\psi}=\frac{1}{2}n^{2}$ and
$\pm$ admissible $\psi$ such that $Q=Q(\partial _{\psi})$.

\vspace{3mm}

\noindent {\it Observation 2.2.} From theorem 3.1 in   [12] -- we infer that
generalised differential operators  give rise to   subalgebras  $\sum _{Q}$
of linear maps (plus zero  map of course) commuting with  a given
generalised difference-tial operator $Q$. The intersection of two
different algebras $\sum _{Q_{1}}$ and  $\sum _{Q_{2}}$  is just zero map
added.

The importance of the above Observation 2.2 as well as the
definition below may be further fully appreciated in the context of  the
Theorem 2.1 and the Proposition 3.1 to come.

\vspace{3mm}

\noindent{\it Definition 2.6.}  Let $\{p_{n} \}_{n \geq 0}$  be  the
normal polynomial sequence {\rm [12]} i.e. $p_{0}(x)=1$ and  $p_{n}(0)=0$;
$n \geq 1$. Then we call it the {\it$\psi$-basic sequence} of the generalised
difference-tial operator $Q$ if in addition $Qp_{n}=n_{\psi}p_{n-1}$. We
shall then call $Q$ in short the $\psi$-difference-tial operator. If $\{
s_{n}\}_{n \geq 0}$ is such a polynomial sequence that
$Qs_{n}=n_{\psi}s_{n-1}$, where $Q$ is the $\psi$-difference-tial operator
then we call $\{ s_{n} \}_{n \geq 0}$ the generalised Sheffer
$Q$-$\psi$-sequence {\rm [13]}. Parallely we define a linear map
$\hat{x}_{Q}:P\to P$ such that $\hat{x}_{Q}p_{n}
=\frac{(n+1)}{(n+1)_{\psi}}p_{n+1}$, $n \geq 0$. We call $\hat{x}_{Q}$ the
operator dual to $Q$. For $Q=Q(\partial _{\psi})=\partial _{\psi}$ we write
for short $\hat{x}_{Q(\psi)}\equiv \hat{x}_{\partial _{\psi}}\equiv
\hat{x}_{\psi}$.

\vspace{3mm}

Let us observe that  $[Q,\hat{x}_{Q}]={\rm id}$  therefore  triples of
operators $\{Q,\hat{x}_{Q},{\rm id}\}$ provide us with a continuous family
of generators of GHW in -- as we call it -- $Q$-representation of
Graves-Heisenberg-Weyl algebra. In the following we shall restrict mostly
to special case of generalised differential operators $Q$, namely to
those $Q=Q(\partial _{\psi})$  which are $\partial _{\psi}$-shift invariant
operators [7, 8, 16]. Let us then start with  appropriate Leibnitz
$\psi$-rules for corresponding $\psi$-derivatives recalling at first  that
{\it admissible sequence} is such a sequence $\psi = \{\psi(q)\}_{n \geq 0}$
that $\psi _{n}(q) \neq 0$; $n \geq 0$ and $\psi _{0}(q)=1$. Then we have
the following obvious observation. If $\psi =\{\psi _{n}(q)\}_{n \geq 0}$
and $\varphi=\{\varphi _{n}(q)\}_{n \geq 0}$ are two admissible sequences
then $[\partial _{\psi}, \partial _{\varphi}]=0$ if $\psi = \varphi$. For
\linebreak any $\psi$-derivative an appropriate Leibnitz $\psi$-rules come
true. Namely it is easy to see that the following Leibnitz $\psi$-rules
hold for any formal series $f$ and $g$ (see Remark 2.1):

\noindent \parbox[c]{4.5cm}{for $\partial _{q}$:} $\partial _{q} (f
\cdot g)=(\partial _{q}f) \cdot g+(\hat{Q}f) \cdot (\partial _{q}g);$

\noindent \parbox[c]{4.5cm}{for $\partial _{R}=R(q\hat{Q}) \partial
_{0}$:} $\partial _{R} (f \bullet g)(z)=R(q\hat{Q}) \{ (\partial _{0}f)(z)
\bullet g(z)+f(0)(\partial _{0}g)(z)\}$

\noindent where -- note -- \quad
$R(q\hat{Q})x^{n-1}=n_{R}x^{n-1};\ (n_{\psi}=n_{R}=n_{R}(q)=R(q^{n}))$
\hspace{1mm} and finally

\noindent \parbox[c]{4.5cm}{for  $\partial
_{\psi}=\hat{n}_{\psi}\partial_{0}$:} $\partial _{\psi}(f \bullet
g)(z)=\hat{n}_{\psi}\{(\partial _{0}f)(z) \bullet g(z)+f(0)(\partial
_{0}g)(z)\}$

\noindent where $\hat{n}_{\psi}x^{n-1}=n_{\psi}x^{n-1}; \quad n \geq 1$.

\vspace{3mm}

\noindent {\it Example 2.2.} [13, 16]  Let  $Q(\partial _{\psi})=D\hat{x}D$,
where $\hat{x}f(x) = x f(x)$ and $D={d}/{dx}$. Then  $\psi =\left\{
[(n^{2})!]^{-1} \right\}_{n \geq 0}$ and $Q=\partial _{\psi}$. Let
$Q(\partial _{\psi})=R(q \hat {Q}) \partial _{0} \equiv \partial _{R}$.
Then $\psi = \left\{\left[ R(q^{n})! \right]^{-1}  \right\}_{n \geq 0}$ and
$Q=\partial _{\psi} \equiv \partial _{R}$. Here $R(z)$ is any formal
Laurent series; $\hat{Q}f(x)=f(qx)$ and $n_{\psi}=R(q^{n})$. $\partial
_{0}$ is $q = 0$ Jackson derivative which as a matter of fact -- being a
difference operator is the differential operator of infinite order at the
same time:
\begin{equation}
\partial _{0}=\sum_{n=1}^{\infty}
(-1)^{n+1}\frac{x^{n-1}}{n!}\frac{d^{n}}{dx^{n}}.
\end{equation}

\noindent Naturally   with the choice $\psi _{n}(q)=[R(q^{n})!]^{-1}$  and
$R(x)=\frac{1-x}{1-q}$ the $\psi$-derivative $\partial _{\psi}$ becomes the
Jackson's derivative [16, 15] $\partial _{q}$: $(\partial _{q}
\varphi)(x)=\frac{1-q\hat{Q}}{(1-q)} \partial_{0}\varphi (x)$. The
form equivalent to (2.6) form of Bernoulli-Taylor expansion one may find in
{\it Acta Eruditorum} from November 1694 [30]  under the name {\it ``series
universalissima"}. (Taylor's expansion was presented in his ``Methodus
incrementorum directa et inversa" in 1715  -- edited in London.) The
example of  $\psi$-derivative:  $Q(\partial _{\psi})=R(q\hat{Q})\partial
_{0} \equiv  \partial _{R}$ i.e. $\psi = \left\{[R(q^{n})!]^{-1}
\right\}_{n \geq 0}$ one may find in [31] where an advanced theory of general
quantum coherent states is being developed. The operator
$R(q\hat{Q})\partial _{0}$ is not recognised in [31] as an example of
$\psi$-derivative. Another very important example of $\psi$-derivative one
may find in [32] where all sequences of binomial type with persistent roots
are classified. There $Q=\partial _{\psi}=2[2D \hat{x}D-D]$ defines -- what
is called by the authors of [32] -- the hyperbolic umbral calculus. Here
$\psi = \left\{ [2n(2n-1)!]^{-1} \right\}_{n \geq 0}$.

Naturally,  choosing any invertible formal series  $S=S(\partial _{\psi})$
in powers of $\partial _{\psi}$ one arrives at infinity of examples of
Appel $\psi$-sequences $s_{n}(x)=S^{-1}(\partial _{\psi})x^{n}$; $n \geq
0$. We are now in a position to define the basic objects of ``$\psi$-umbral
calculus" [16, 15].

\hspace{3mm}

\noindent {\it Definition 2.7.} Let $Q(\partial _{\psi}):P\to P$;  the
linear operator $Q(\partial _{\psi})$ is a {\it $\partial _{\psi}$-delta
operator} if

a) $Q(\partial _{\psi})$  is $\partial _{\psi}$-shift invariant;
\hspace{3mm} b) $Q(\partial _{\psi})(id)={\rm const}\neq 0$.

The strictly related notion is that of the  $\partial _{\psi}$-basic
polynomial sequence:

\vspace{3mm}

\noindent {\it Definition 2.8.} Let $Q(\partial _{\psi}):P\to P$;  be
the  $\partial _{\psi}$-delta operator.  A polynomial sequence
$\{p_{n}\}_{n \geq 0}$; $\deg p_{n}= n$ such that: {\rm 1)}
$p_{0}(x)=1;$\hspace{3mm} {\rm 2)} $p_{n}(0)=0;\ n>0;$\hspace{3mm} {\rm 3)}
$Q(\partial _{\psi})p_{n}=n_{\psi}p_{n-1}$ is called the {\it $\partial
_{\psi}$-basic polynomial sequence} of the $\partial _{\psi}$-delta operator
$Q(\partial _{\psi})$.

\vspace{3mm}

\noindent {\it Note 2.1.}    Let  $\Phi (x;\lambda)=\sum_{n \geq
0}\frac{\lambda ^{n}}{n_{\psi}!}p_{n}(x)$ denotes the $\psi$-exponential
generating function of the $\partial _{\psi}$-basic polynomial sequence
$\{p_{n} \}_{n \geq 0}$ of the $\partial _{\psi}$-delta operator $\hat{Q}
\equiv Q(\partial \psi)$ and let $\Phi(0;\lambda)=1$. Then $\hat{Q}
\Phi(x;\lambda)=\lambda \Phi (x;\lambda)$ and $\Phi$ is the unique solution
of this eigenvalue problem. In view of the Observation 2.2 we affirm that
there exists such an admissible sequence $\varphi$ that
$\Phi(x;\lambda)=\exp _{\varphi} [\lambda x]$.

The
notation and naming established by  Definitions 2.7 and 2.8  serve the
target  to preserve and to broaden simplicity of Rota's finite operator
calculus also in its extended ``$\psi$-umbral calculus" case [16, 15]. As a
matter of illustration of such notation efficiency  let us quote after
[16]  the important Theorem 2.1  which might be proved using the fact that
$\forall\ Q(\partial _{\psi}) \exists !$ invertible  $S \in \sum _{\psi}$
 such that $Q(\partial _{\psi})=\partial _{\psi}S$. (For Theorem 2.1
see also Theorem 4.3 in [12] -- which holds for operators introduced by
the definitions 2.5 and 2.6.)

\vspace{3mm}

\noindent {\bf Theorem 2.1.}  (Lagrange  and Rodrigues  $\psi$-formulas)

{\it Let $\{p_{n}(x)\}_{n=0}^{\infty}$   be  $\partial _{\psi}$-basic polynomial sequence
of the $\partial _{\psi}$-delta operator $Q(\partial _{\psi})$. Let
$Q(\partial _{\psi})=\partial _{\psi}S$. Then for $n>0$}:

(1) $p_{n}(x) =  Q(\partial _{\psi})'S^{-n-1}x^{n}$;

(2) $p_{n}(x)=S^{-n}x^{n}-\frac{n_{\psi}}{n}(S^{-n})'x^{n-1}$;

(3) $p_{n}(x)=\frac{n_{\psi}}{n}\hat{x}_{\psi}S^{-n}x^{n-1}$;

(4) $p_{n}(x)=\frac{n_{\psi}}{n}\hat{x}_{\psi}\left( Q(\partial _{\psi})'
\right)^{-1}p_{n-1}(x)$\hspace{3mm} (Rodrigues  $\psi$-formula)

\noindent For the proof  one  uses  typical properties of  the Pincherle
$\psi$-derivative as well as  $\hat{x}_{\psi}$   operator  [16]   defined
below.

\vspace{3mm}

\noindent {\it Definition 2.9.}  (compare with  (17) in  [7])   The
{\it Pincherle $\psi$-derivative} i.e. the  linear map    $': \sum_{\psi} \to \sum
_{\psi};$ is given by  $T'=T\hat{x}_{\psi}-\hat{x}_{\psi}T_{\partial
_{\psi}}\equiv [T,\hat{x}_{\psi}]$ where the linear map $\hat{x}_{\psi}: P
\to P;$ is defined in the basis $\{ x^{n}\}_{n \geq 0}$ as follows

\[
\hat{x}_{\psi}x^{n}=\frac{\psi_{n+1}(q)(n+1)}{\psi_{n}(q)}\
x^{n+1}=\frac{(n+1)}{(n+1)_{\psi}}\ x^{n+1};\quad n \geq 0. \]

\vspace{3mm}

\noindent {\it Observation 2.3.}  The triples  $\{\partial _{\psi},
\hat{x}_{\psi},{\rm id}\}$ for any admissible $\psi$  constitute  set of
generators of the $\psi$-labelled representations of
Graves-Heisenberg-Weyl (GHW) algebra [33--35]. Namely, as easily seen
$[\partial _{\psi}, \hat{x}_{\psi}]={\rm id}$.

\vspace{3mm}

\noindent {\it Observation 2.4.} In view of the
Observation 2.3 the general Leibnitz rule in \linebreak $\psi$-representation of
Graves-Heisenberg-Weyl algebra may be written (compare with  2.2.2
Proposition in [35])  as follows

\begin{equation}
\partial _{\psi}^{n}\hat{x}_{\psi}^{m}=\sum _{k \geq 0}\left(
\begin{array}{c} n\\k \end{array} \right)\left(
\begin{array}{c} m\\k \end{array} \right) k!\ \hat{x}_{\psi}^{m-k}\ \partial
_{\psi}^{n-k}. \end{equation}

\noindent One derives the above $\psi$-Leibnitz rule from
$\psi$-Heisenberg-Weyl exponential commutation rules exactly the same way
as in $\{D,\hat{x}, {\rm id} \}$ GHW representation -- (compare with  2.2.1
Proposition in [35]). $\psi$-Heisenberg-Weyl exponential commutation
relations read:

\begin{equation}
\exp\{t\partial
_{\psi}\}\exp\{a\hat{x}_{\psi}\}=\exp\{at\}\exp\{a\hat{x}_{\psi}\}\exp \{t
\partial _{\psi}\}. \end{equation}

In the sequel we shall discover another, related to  (2.6)  and
heuristically appealing  form of a  $\partial _{\psi}$-difference-ization
rule for a specific  new  $*_{\psi}$   product of  functions or formal series.  To
this end let us introduce this pertinent  $\psi$-multiplication  $*_{\psi}$   of
functions or formal series as specified below.

\vspace{3mm}
\noindent {\it Notation 2.2.}

\noindent
$x*_{\psi}x^{n}=\hat{x}_{\psi}(x^{n})=\frac{(n+1)}{(n+1)_{\psi}}x^{n+1};
\quad n \geq 0\quad {\rm hence}\quad x*_{\psi}1=1_{\psi}xTx$

\noindent
$x^{n}*_{\psi}x=\hat{x}_{\psi}^{n}(x)=\frac{(n+1)!}{(n+1)_{\psi}!}x^{n+1};
\quad n \geq 0\quad {\rm hence}\quad 1*_{\psi}x=1_{\psi}xTx$ therefore
$x*_{\psi}\alpha$

\noindent $1=\alpha 1*_{\psi}x=x*_{\psi}\alpha=\alpha *_{\psi}x=\alpha
1_{\psi}x\quad \forall x, \ \alpha \in \otimes;\quad
f(x)*_{\psi}x^{n}=f(\hat{x}_{\psi})x^{n}$.

For $k \neq n,$ $x^{n}*_{\psi}x^{k}\neq x^{k}*_{\psi}x^{n}$ as well as
$x^{n}*_{\psi}x^{k}\neq x^{n+k}$ -- in general i.e. for arbitrary
admissible $\psi$; compare this with $(x +_{\psi} a)^{n}\neq (x
+_{\psi}a)^{n-1}(x +_{\psi} a)$.

In order to facilitate in the future
formulation of observations accounted for  on the basis of $\psi$-calculus
representation of GHW algebra we shall use what follows.

\vspace{3mm}

\noindent {\it Definition 2.10.}  With  Notation 2.2.  adapted let
us define the $*_{\psi}$  {\it powers} of   $x$  according to

\[
x^{n^{*}_{\psi}}\equiv x*_{\psi}x^{(n-1)^{*}_{\psi}}=\hat{x}_{\psi}\left(
x^{(n-1)^{*}_{\psi}}\right)=x*_{\psi}x*_{\psi}...*_{\psi}x=\frac{n !}
{n_{\psi} !}x^{n}; \quad n \geq 0.\]

\noindent  Note that
$x^{n^{*}_{\psi}}*_{\psi}x^{k^{*}_{\psi}}=\frac{n!}{n_{\psi}!}x^{(n+k)^{*}_
{\psi}} \neq
x^{k^{*}_{\psi}}*_{\psi}x^{n^{*}_{\psi}}=\frac{k!}{k_{\psi}!}x^{(n+k)^{*}
_{\psi}}$ for $k \neq n$ and \linebreak $x^{0^{*}_{\psi}}=1$. This
noncommutative $\psi$-product $*_{\psi}$ is devised so as to ensure the
observations below.

\vspace{3mm}

\noindent {\it Observation 2.5.}

\noindent {\it a)}\hspace{0.5cm} $\partial
_{\psi}x^{n^{*}_{\psi}}=nx^{(n-1)^{*}_{\psi}};\quad n \geq 0$

\noindent {\it b)}\hspace{0.5cm} $\exp _{\psi}[\alpha x]\equiv \exp \{ \alpha
\hat{x}_{\psi} \} {\bf 1}$

\noindent {\it c)}\hspace{0.5cm} $\exp[\alpha x]*_{\psi}\left\{ \exp
_{\psi}\{\beta \hat{x}_{\psi} \}{\bf 1} \right\} = \exp _{\psi}\left\{
[\alpha + \beta] \hat{x}_{\psi} \right\} {\bf 1}$

\noindent {\it d)}\hspace{0.5cm} $\partial _{\psi}\left(
x^{k}*_{\psi}x^{n^{*}_{\psi}} \right) =(D
x^{k})*_{\psi}x^{n^{*}_{\psi}}+x^{k}*_{\psi}\left( \partial
_{\psi}x^{n^{*}_{\psi}} \right) $

\noindent {\it e)}\hspace{0.5cm}  $\partial _{\psi}\left( f*_{\psi}g
\right)=(Df)*_{\psi}g+f*_{\psi}(\partial _{\psi}g); \quad f,g$ -- {\it
formal series}

\noindent {\it f)}\hspace{0.5cm} $f(\hat{x}_{\psi})g(\hat{x}_{\psi}){\bf
1}=f(x)*_{\psi}\tilde{g}(x);\quad \tilde{g}(x)=g(\hat{x}_{\psi}){\bf 1}.$

Now the consequences of Leibnitz rule  e)  for  difference-ization of the
product  are easily feasible.  For example the Poisson $\psi$-process
distribution   $\pi_{m}(x)=\frac{1}{N(\lambda,x)}p_{m}(x)$; $\sum_{m
\geq0}\pi_{m}(x)=1$ is determined by

\begin{equation}
p_{m}(x)=\frac{(\lambda x)^{m}}{m!}*_{\psi}\exp _{\psi}\left[ -\lambda x
\right] \end{equation}

\noindent which is the (up to a factor) unique solution of its
corresponding  $\partial _{\psi}$-difference equations system

\begin{equation}
\partial _{\psi}p_{m}(x)+\lambda p_{m}(x)=\lambda p_{m-1}(x)m>0; \quad
\partial _{\psi}p_{0}(x)=-\lambda p_{0}(x) \end{equation}

\noindent Naturally $N(\lambda , x)=\exp [\lambda x]*_{\psi}\exp
_{\psi}[-\lambda x]$.

As announced -- the  rules  of  $\psi$-product   $*_{\psi}$   are accounted
for -- as a matter of fact -- on the basis of  $\psi$-calculus
representation of GHW algebra.  Indeed, it is enough to consult
Observation 2.3 and to introduce  $\psi$-Pincherle derivation
$\hat{\partial}_{\psi}$ of series in powers of  the symbol $\hat{x}_{\psi}$
as below. Then the correspondence between generic relative  formulas turns
out evident.

\vspace{3mm}

\noindent {\it Observation 2.6.}  Let  $\hat{\partial}_{\psi}\equiv
\frac{\partial}{\partial \hat{x}_{\psi}}$ be defined according to
$\hat{\partial}_{\psi}f(\hat{x}_{\psi})=[\partial _{\psi},
f(\hat{x}_{\psi})]$. Then

\[
\hat{\partial}_{\psi}\hat{x}^{n}_{\psi}=n\hat{x}^{n-1}_{\psi};\quad n \geq 0 \quad
and\quad \hat{\partial} _{\psi}\hat{x}^{n}_{\psi}{\bf 1}=\partial
_{\psi}x^{n^{*}_{\psi}} \quad hence \quad \left[
\hat{\partial}_{\psi}f\left( \hat{x}_{\psi} \right) \right]{\bf 1}=\partial
_{\psi}f(x); \]

\noindent {\it where $f$ is a formal series
in powers $\hat{x}_{\psi}$ of or equivalently in $*_{\psi}$ powers of $x$}.

As an  example of application  note how   the solution of   (2.10)  is
obtained from the obvious solution $\pi _{m}\left( \hat{x}_{\psi} \right)$
of the $\hat{\partial}_{\psi}$-Pincherle differential equation  (2.11)
formulated within G-H-W algebra generated by $\{\partial _{\psi},
\hat{x}_{\psi}, {\rm id} \}$

\begin{equation}
\hat{\partial} _{\psi} \pi _{m}\left( \hat{x}_{\psi} \right) +\lambda \pi
_{m}(\hat{x}_{\psi})=\lambda \pi
_{m-1}(\hat{x}_{\psi})m>0; \quad \partial
_{\psi}\pi_{0}(\hat{x}_{\psi})=-\lambda \pi _{0}(\hat{x}_{\psi})
\end{equation}

\noindent Namely:   due to Observation 2.5 f)  $p_{m}(x) = \pi
_{m}(\hat{x}_{\psi} ){\bf 1}$, where

\begin{equation}
\pi _{m}(\hat{x}_{\psi})= \frac{\left( \lambda \hat{x}_{\psi}
\right)^{m}}{m!}\exp_{\psi}\left[ -\lambda \hat{x}_{\psi} \right].
\end{equation}


\subsection{The general picture of  $End(P)$}
\setcounter{equation}{0}

In the following we shall consider the algebra  $End(P)$ of linear
operators on  linear space $P$ of polynomials  $P = {\bf F}[x]$ over the
field ${\bf F}$ of characteristic zero.

Before we pass to
characterisations of Sheffer $\psi$-polynomials and related propositions we
shall draw  an overview picture of the situation with underlying the
possibility to develop umbral calculus for  {\it any} polynomial sequences
$\{ p_{n} \}^{\infty}_{0}$ instead of those of traditional binomial type
only. At first note that in 1901 it was proved [36]  that every linear
operator mapping $P$ into  $P$ may be represented as infinite series in
operators $\hat{x}$  and $D$. In 1986 the authors of  [37] supplied an
explicit expression for such series in most general case of polynomials in
one variable (for many variables see [38]). Thus according to Proposition
1 from [37] one has:

\vspace{3mm}

\noindent {\bf Proposition 3.1.}  {\it Let $\hat{Q}$  be a linear operator
that reduces by one the degree of  each polynomial. Let $\{ q_{n}(\hat{x})
\}_{n \geq 0}$ be an arbitrary sequence of polynomials in the operator
$\hat{x}$. Then $\hat{T}=\sum _{n \geq 0}q_{n}(\hat{x}) \hat{Q}^{n}$
defines a linear operator that maps polynomials into polynomials.
Conversely, if $\hat{T}$ is linear operator that maps polynomials into
polynomials then there exists a unique expansion of the form
$\hat{T}=\sum_{n \geq 0}q_{n}(\hat{x})\hat{Q}^{n}$}.

\vspace{2mm}

\noindent It is also a rather matter of an easy  exercise (see Note 2.1.) to prove [37]:

\vspace{3mm}

\noindent {\bf Proposition 3.2.}  {\it Let  $\hat{Q}$ be a linear operator that
reduces by one the degree of each polynomial. Let $\{ q_{n}(\hat{x})
\}_{n \geq 0}$   be an arbitrary
sequence of polynomials in the operator   $\hat{x}$. Let  a linear operator
that maps polynomials into polynomials be given by   $\hat{T}=\sum_{n \geq
0}q_{n}(\hat{x})\hat{Q}^{n}$. Let $P(x;\lambda)=\sum_{n \geq
0}q_{n}(x)\lambda ^{n}$ be the indicator of $\hat{T}$. Then there exists a
unique formal series $\Phi(x;\lambda);\ \Phi(0;\lambda)=1$ such that
$\hat{Q}\Phi(x;\lambda)=\lambda \Phi(x; \lambda)$ and
$P(x;\lambda)=\Phi(x;\lambda)^{-1}\hat{T}\Phi(x;\lambda)$}.

\vspace{3mm}

\noindent {\it Example 3.1.} \hspace{5mm} $\partial _{\psi}\exp
_{\psi}\{\lambda x \}=\lambda \exp_{\psi}\{\lambda x \}; \quad
\exp_{\psi}[x]_{x=0}=1$. \hfill (*)

\vspace{2mm}

\noindent Hence for indicator of   $\hat{T}$; $\hat{T}=\sum_{n \geq 0}
q_{n}(\hat{x})\partial ^{n}_{\psi}$ we have

\vspace{2mm}

\hspace{3cm}$P(x;\lambda)=\left[\exp_{\psi}\{\lambda x
\}\right]^{-1}\hat{T}\exp _{\psi}\{\lambda x\}.$\hfill(**)

\vspace{2mm}

\noindent After choosing  $\psi _{n}(q)=\left[ n_{q}! \right]^{-1}$  we get
$\exp _{\psi}\{x\} = \exp _{q} \{x\}$. In this connection note that
$\exp_{0}(x)=\frac{1}{1-x}$ and $\exp(x)$ are mutual limit deformations.
Therefore corresponding specifications  of (*) such as $\exp_{0}(\lambda
x)=\frac{1}{1-\lambda x}$ or $\exp(\lambda x)$ lead to corresponding
specifications of (**) for divided difference operator $\partial _{0}$ and
$D$ operator including special cases from [37].

Let us now introduce [16, 15]  an  important operator $\hat{x}_{Q(\partial
_{\psi})}$ dual to $Q(\partial_{\psi})$.

\vspace{3mm}

\noindent {\it Definition 3.1.}  Let $\{p_{n}\}_{n \geq 0}$  be  the
$\partial _{\psi}$-basic polynomial sequence of the  $\partial
_{\psi}$-delta operator $Q(\partial _{\psi})$.  A linear map
$\hat{x}_{Q(\partial _{\psi})}:P \to P$; $\hat{x}_{Q(\partial
_{\psi})}p_{n}=\frac{(n+1)}{(n+1)_{\psi}}p_{n+1}$; $n \geq 0$
is called the operator {\it dual} to  $Q(\partial _{\psi})$. (Note:  for  $Q =
{\rm id}$ we have $\hat{x}_{Q(\partial _{\psi})}\equiv \hat{x}_{\partial
_{\psi}}\equiv \hat{x}_{\psi}$.)

\vspace{3mm}

\noindent {\it Comment 3.1.}  Dual in the above sense corresponds to
adjoint in $\psi$-umbral calculus language of  linear functionals' umbral
algebra   (compare with  Proposition 1.1.21 in [18]).

It is obvious that the following holds:

\vspace{3mm}

\noindent {\bf Proposition 3.3.} {\it  Let $\{q_{n}(\hat{x}_{Q(\partial
_{\psi})}) \}_{n \geq 0}$ be an arbitrary sequence of polynomials in the
operator $\hat{x}_{Q(\partial _{\psi})}$. Then $T=\sum _{n \geq
0}q_{n}(\hat{x}_{Q(\partial _{\psi})})Q(\partial _{\psi})^{n}$ defines a
linear operator that maps polynomials into polynomials. Conversely, if $T$
is linear operator that maps polynomials into polynomials then there exists
a unique expansion of the form}
\begin{equation}
T=\sum _{n \geq 0}q_{n}\left( \hat{x}_{Q(\partial _{\psi})}\right)Q(\partial
_{\psi} )^{n} \end{equation}

\noindent {\it Comment 3.2.}  The pair $Q(\partial _{\psi})$,
$\hat{x}_{Q(\partial _{\psi})}$ of dual operators is expected to play a
role in the description of quantum-like processes apart from the $q$-case
now vastly exploited [19, 15, 16]. Naturally the Proposition 3.2 for
$Q(\partial _{\psi})$ and $\hat{x}_{Q(\partial _{\psi})}$ dual operators is
also valid.

\vspace{3mm}

\noindent {\bf Summing up:}  (see Definitions 2.5  and  2.6) We have the
following picture for  $End (P)$ -- the algebra of all linear operators
acting on the algebra  $P$ of polynomials: $Q(P ) \subset End (P)$ and of
course $Q(P ) \neq End (P)$    where $Q(P)$  is the linear space of  all
$\psi$-difference-tial operators $Q$  (including   $\partial _{\psi}$-delta
operators $Q(\partial _{\psi})$ and plus zero map of course).

$Q(P)$ breaks
up into sum [12] of subsets $\sum _{L} \cap\  Q(P)^{*}$; $L \in Q(P )$
according to commutativity of these generalised differential operators $L$
where  $\sum_{L}=\left\{T \in End(P\right.)$; $\left.[T,L]=0;\ L\in Q(P)
\right\}$. Any one dimensional subspace of $End(P)$ spanned by arbitrary
operator $\hat{r} \in End(P)$ raising by one the degree of all polynomials
has the empty intersection with $Q(P)^{*} = Q(P) / \{0\}$. Naturally to
each $\sum _{\psi}$ subalgebra i.e. to each $Q$ operator there corresponds
its dual operator $\hat{x}_{Q}$; $\hat{x}_{Q} \notin \sum_{Q}$ and both $Q$
and $\hat{x}_{Q}$ operators are sufficient to build up the whole algebra
$End(P)$ according to unique representation given by (3.1). With every such
dual pair  $\{Q, \hat{x}_{Q}\}$ including $Q(\partial _{\psi})$ and
$\hat{x}_{Q(\partial _{\psi})}$ which include $Q = {\rm id.}$ case i.e.
$\partial _{\psi}$ and $\hat{x}_{\partial _{\psi}}\equiv \hat{x}_{\psi}$ we
have for any admissible $\psi$ the general picture statement.

\vspace{3mm}

\noindent {\bf General statement:}
\[
End(P) = \left[\left\{\partial _{\psi}, \hat{x}_{\psi} \right\}\right] =
\left[\left\{Q(\partial _{\psi}), \hat{x}_{Q(\partial _{\psi})}
\right\}\right]=\left[\left\{Q, \hat{x}_{Q} \right\}\right] \]

\noindent i.e.  the algebra   $End(P)$ is generated by any dual pair $\{Q,
\hat{x}_{Q} \}$ and in particular   $End (P)$  is generated by
$\{Q(\partial _{\psi}), \hat{x}_{Q(\partial _{\psi})} \}$ including any
dual pair $\{\partial _{\psi}, \hat{x}_{\psi}\}$ determined by any choice of
admissible sequence $\psi$.

As a matter of fact and in another words: we have bijective
correspondences between  different commutation classes  of operators from
$End (P)$; different abelian \linebreak subalgebras  $\sum _{Q}$; distinct
$\psi$-representations of GHW algebras; different $\psi$-representa-tions of
the reduced incidence algebra $R(L(S))$ -- isomorphic to the algebra
$\Phi_{\psi}$ of $\psi$-exponential formal power series [16]   and finally
-- distinct $Q$-umbral or $\psi$-umbral calculi [7, 8, 12, 13, 32, 16].
(Recall: $R(L(S))$ is the reduced incidence algebra of $L(S)$ where
$L(S)=\{A;\ A\subset S;\ |A|<\infty\};\ S$ is countable and $(L(S);\
\subseteq)$ is partially ordered set ordered by inclusion  [3, 16]).

This is the way the Rota's  devise has been carried into effect. The devise
{\it ``much is the iteration of the few"} [3] -- much of the properties of
literally {\it all} polynomial sequences  -- as well as GHW algebra
representations --  is the application of  few basic principles of the
$Q$-umbral or in particular   $\psi$-umbral difference operator calculus
``$\psi$-{\it integration}"  included.

\vspace{3mm}

\noindent $\psi$-{\bf Integration Remark:} Recall: $\partial
_{0}x^{n}=x^{n-1}$. $\partial _{0}$ is identical with divided difference
operator. $\partial _{0}$ is identical with $\partial _{\psi}$ for
$\psi =\{\psi(q)_{n} \}_{n \geq 0}$; $\psi(q)_{n}=1$; $n \geq 0$. Let
$\hat{Q}f(x) = f(qx)$. Recall also that to the ``$\partial_{Q}$
difference-ization" there corresponds also $q$-integration [24--26] which is
a right inverse operation to ``$q$-difference-ization". Namely
\begin{equation}
F(z):\equiv \left( \int _{q} \varphi \right) (z):=\left( 1-q \right)z\sum
_{k=0}^{\infty}\varphi \left( q^{k}z \right)q^{k}
 \end{equation}

\noindent i.e.
\begin{equation}
F(z)\equiv \left( \int _{q}\varphi \right)(z)=(1-q)z \left( \sum
_{k=0}^{\infty}q^{k}\hat{Q}^{k}\varphi \right)(z)=\left( (1-q)z
\frac{1}{1-q\hat{Q}}\varphi \right)(z)
\end{equation}

\noindent Of course
\begin{equation}
\partial _{q} \circ \int _{q} =id
\end{equation}

\noindent as
\begin{equation}
\frac{1-q\hat{Q}}{1-q}\partial _{0} \left( (1-q)\hat{z}
\frac{1}{1-q\hat{Q}} \right)={\rm id}. \end{equation}

\noindent Naturally (3.5) might serve to define a right inverse operation
to ``$q$-difference-ization"  $(\partial
_{q}\varphi)(x)=\frac{1-q\hat{Q}}{1-q}\partial_{0}\varphi(x)$ and
consequently the ``$q$-integration" as represented by (3.2) and (3.3). As
it is well known the definite $q$-integral is a numerical approximation of
the definite integral obtained in the  $q \to 1$ limit.  Following the
$q$-case example we introduce now an $R$-integration   (consult Remark 2.1).
\begin{equation}
\int _{R} x^{n}=\left( \hat{x}\frac{1}{R(q\hat{Q})} \right)
x^{n}=\frac{1}{R(q^{n+1})}x^{n+1}; \quad n \geq 0 \end{equation}

\noindent Of course $\partial _{R} \circ \int _{R}={\rm id}$ as
\begin{equation}
R(q \hat{Q})\partial _{0} \left( \hat{x}\frac{1}{R(q\hat{Q})} \right)={\rm id}.
\end{equation}

Let us  then finally introduce the analogous representation  for
$\partial _{\psi}$ difference-ization
\begin{equation}
\partial_{\psi}=\hat{n}_{\psi}\partial_{0}; \quad
\hat{n}_{\psi}x^{n-1}=n_{\psi}x^{n-1};\quad n\geq 1.
\end{equation}

\vspace{2mm}

\noindent Then
\begin{equation}
\int_{\psi}x^{n}=\left(\hat{x}\frac{1}{\hat{n}_{\psi}}
\right)x^{n}=\frac{1}{(n+1)_{\psi}}x^{n+1}; \quad n \geq 0 \end{equation}

\vspace{3mm}

\noindent {\bf The Section Closing Remark:}  The picture that emerges
discloses the fact that -- any $\psi$-representation of finite operator
calculus or equivalently -- any \linebreak $\psi$-representation  of GHW algebra makes up
an example of the algebraization of the analysis -- naturally when
constrained to the algebra of polynomials. We did  restricted  all our
considerations to the algebra $P$ of polynomials or formal series.
Therefore the distinction in-between difference and differentiation
operators disappears. All linear operators  on  $P$  are both difference and
differentiation operators if the degree of differentiation or difference
operator is unlimited.  For example
\begin{eqnarray}
\frac{d}{dx}=\sum_{k \geq 1}\frac{d_{k}}{k!} \Delta _{k}& {\rm where}
& d_{k}=\left[ \frac{d}{dx} x^{\underline{k}}
\right]_{x=0}=(-1)^{k-1}(k-1)!\nonumber\\
&& \nonumber \\
& {\rm or}& \Delta =\sum _{n \geq 1}
\frac{\delta _{n}}{n!}\frac{d^{n}}{dx^{n}} \nonumber
\end{eqnarray}

\noindent where  $\delta _{n}=[\Delta x^{n}]_{x=0}=1$. Thus the difference
and differential operators and equations are treated on the same footing.


\subsection{Characterisations of Sheffer $\psi$-polynomials and
related \\ propositions}
\setcounter{equation}{0}

Let us now pass to characterisations of Sheffer $\psi$-polynomials and
related propositions.

\vspace{2mm}

\noindent {\it Definition 4.1.}   A polynomial sequence  $\{s_{n}(x)
\}_{n=0}^{\infty}$ is called the {\it Sheffer   $\psi$-sequence} of  the
$\partial _{\psi}$-delta operator $Q(\partial _{\psi})$ if (1)
$s_{0}(x)=c\neq 0;$ (2) $Q(\partial _{\psi})s_{n}(x)=n_{\psi}s_{n-1}(x)$.
The following characterisation of Sheffer $\psi$-sequence of the
$\partial _{\psi}$-delta operator $Q(\partial _{\psi})$ relates it to the
unique pair of a $\partial _{\psi}$-basic sequence of $Q(\partial _{\psi})$
and the corresponding $\partial _{\psi}$-shift invariant operator $S$ [16].

\vspace{2mm}

\noindent {\bf Proposition 4.1.} {\it Let $Q(\partial _{\psi})$  be a
$\partial _{\psi}$-delta operator with $\partial _{\psi}$-basic polynomial
sequence $\{q_{n}(x)\}^{\infty}_{n=0}$. Then $\{s_{n}(x)\}^{\infty}_{n=0}$
is a sequence of Sheffer $\psi$-polynomials of $Q(\partial _{\psi})$ if
there exists an invertible $\partial _{\psi}$-shift invariant operator
$S$ such  that} $s_{n}(x)=S^{-1}q_{n}(x)$. Sheffer   $\psi$-sequence
$\{s_{n}(x)\}^{\infty}_{n=0}$ labelled by $S$ is therefore referred to as
the Sheffer $\psi$-sequence of $Q(\partial_{\psi})$ relative to  $S$. It is not difficult to
find out that $S^{-1}=\sum_{k \geq 0}\frac{s_{k}(0)}{k_{\psi}!}Q(\partial
_{\psi})^{k}$, where $\{s_{n}(x)\}^{\infty}_{n=0}$ is the sequence of
Sheffer $\psi$-sequence $Q(\partial _{\psi})$ relative to $S$ [16].

\vspace{2mm}

\noindent {\it Remark 4.1} The Proposition 4.1  is also valid in the
``$Q$-case" also   i.e.  for  the  $\psi$-difference-tial operators $Q$
and for the generalised Sheffer  $\psi$-polynomials  $\{ s_{n}\}_{n \geq
0}$ (see definitions 2.5 and 2.6). In such a  $Q$-case we shall say that
$Q$-$\psi$-sequence $\{s_{n}(x)\}_{n=0}^{\infty}$ is the generalised
Sheffer $\psi$-sequence of $Q$ relative to $S$.

\vspace{2mm}

With this in mind we arrive at the general conclusion.

\vspace{2mm}

\noindent {\bf Conclusion 4.1.}  The family of generalised Sheffer
$\psi$-sequences $\{s_{n}(x) \}^{\infty}_{n=0}$  corresponding to a fixed
$\psi$-difference-tial operator $Q$ i.e. the family of   Sheffer
$Q$-$\psi$-sequences (see the  definitions 2.5  and 2.6)  is labelled by
the abelian group of all invertible operators $S\in \Sigma _{Q}$. The
families of these generalised Sheffer $\psi$-polynomials are orbits of such
groups.

$\psi$-calculus or $Q$-calculus in parts appears to be almost automatic
extension of the classical operator calculus of Rota-Mullin [16].
Therefore -- as a rule -- we shall omit proofs of propositions and theorems
stated.  Proofs might be either found in [16, 15, 39, 40] or  adopted from
corresponding ones in [3] in most of relevant cases.  Recall again that the
general results of $\psi$-calculus presented here may be extended to
Markowsky $Q$-calculus where  $Q$ stands for a generalised difference
operator i.e. the one lowering the degree of any polynomial by one.

As a matter of fact all statements of standard finite operator calculus of
Rota  are valid also in the case of $\psi$-extension under the almost
automatic replacement \linebreak of  $\{D, \hat{x}, {\rm id} \}$ by their
$\psi$-representation correspondent generators of  GHW  i.e. \linebreak $\{\partial
_{\psi}, \hat{x}_{\psi}, {\rm id} \}$. In most general case of $Q$-umbral
representation these are the \linebreak $\{Q, \hat{x}_{Q}, {\rm id} \}$ triples of
operators to be used. Naturally any specification of admissible $\psi$ --
for example the famous one defining $q$-calculus -- has its own
characteristic properties not pertaining to the standard case of Rota
calculus realisation. Nevertheless the overall picture and system of
statements depending only on  GHW algebra is the same modulo some automatic
replacements  in formulas  quoted in the sequel [16, 39].  Let us then pass
to the presentation of characterisations of Sheffer  $\psi$-sequences  and
related propositions in order to indicate the first main features of
extended finite operator calculus. Apart from Sheffer $\psi$-sequences
characterisation by Proposition 4.1 and Viskov theorems from section 2
-- (see Illustration 2.1 -- formulas (2.1--2.4)) the following
propositions and theorems come true [16, 15, 39]. At first  we restate the
Sheffer  $\psi$-sequences characterisation due to Viskov (see Illustration
2.1.) in our Rota-oriented further notation [16].

\vspace{2mm}

\noindent {\bf Proposition 4.2.} {\it Let $Q$  be a  $\partial
_{\psi}$-delta operator. Let $S$ be an invertible  $\partial _{\psi}$-shift
invariant operator. Let $Q=q(\partial _{\psi})$ and $S=s(\partial
_{\psi})$. Let $q^{-1}(t)$ be the $\psi$-exponential formal power series
inverse to $q(t)$. Then the $\psi$-exponential generating function of
Sheffer $\psi$-sequence $\{ s_{n}(x)_{n=0}^{\infty} \}$ of $Q$ relative to
$S$ is given by}
\[ \sum_{k \geq 0}\frac{s_{k}(x)}{k _{\psi}!}z^{k}=s\left( q^{-1}(z) \right)
\exp _{\psi}\left\{ xq^{-1}(z) \right\}. \]

\noindent Naturally the corresponding extended versions of binomial
and second binomial theorems
hold and here we quote the second one as another  Sheffer
$\psi$-sequences characterisation.

\vspace{2mm}

\noindent {\bf Theorem  4.1.}     (Sheffer  $\psi$-Binomial
Characterisation Theorem)

\noindent {\it Let $Q(\partial _{\psi})$  be the  $\partial _{\psi}$-delta
operator with the $\partial _{\psi}$-basic polynomial sequence   $\{
q_{n}(x) \}^{\infty}_{n=0}$. Then $\{
s_{n}(x) \}^{\infty}_{n=0}$ is the Sheffer $\psi$-sequence of $Q(\partial
_{\psi})$ relative to an invertible $\partial _{\psi}$-shift invariant
operator $S$ if}
\[
s_{n}\left(x+_{\psi}y \right)= \sum _{k \geq 0}\left( \begin{array}{c} n\\k
\end{array} \right)_{\psi} s_{k} (x)q_{n-k}(y). \]

\noindent

There exist  other  characterisations of course and we are not going to
present all of available ones. We  finish our examples of such
characterisation statements with another  proposition characterising
Sheffer   $\psi$-sequences [3, 16].

\vspace{2mm}

\noindent {\bf Proposition 4.3.}   {\it A sequence  $\{s_{n}(x)
\}^{\infty}_{n=0}$  is  the Sheffer $\psi$-sequence of  a  $\partial
_{\psi}$-delta operator $Q(\partial _{\psi})$ with the  $\partial
_{\psi}$-basic polynomial sequence $\{
q_{n}(x) \}^{\infty}_{n=0}$ if there exists such a $\partial
_{\psi}$-delta operator $A$ (not necessarily associated with $\{
s_{n}(x) \}^{\infty}_{n=0}$) and  the
sequence $\{
c_{n}\}^{\infty}_{n=0}$ of constants such that}
\[
As_{n}(x)=\sum _{k \geq 0} \left( \begin{array}{c} n \\k \end{array} \right)
s_{k}(x)c_{n-k}\quad n \geq 0.
\]

We shall present now some propositions announced earlier which are related
to Sheffer  $\psi$-sequences. We have chosen for that the spectral theorem
and few, main umbral propositions [3, 16, 40]. To start with let us note
that the natural inner product may be associated with any
$Q$-$\psi$-sequence $\{ s_{n}(x) \}^{\infty}_{n=0}$ of a
$\psi$-difference-tial operator $Q$ -- relative to $S \in \Sigma _{Q}$. For
that to see define the linear operator $W$:  $s_{n}(x) \to x^{n}$.

\vspace{2mm}

\noindent {\it Definition 4.2.}  Let  $S$   be  an operator the
Sheffer $\psi$-sequence  $\{s_{n}(x) \}_{n=0}^{\infty}$  of generalised
differential operator $Q$ is related to. Let  $W$ be the linear operator
such that  $W:$ $s_{n}(x)\# x^{n}$. We then define the bilinear form
\[
\left ( f(x),g(x) \right)_{Q,S}:=\left[(Wf)(Q)Sg(x)  \right] _{x=0}\quad
f,g \in P. \]

\noindent It is then easy to observe  the important property of this
bilinear form -- now on the reals.

\vspace{2mm}

\noindent {\it Observation 4.1.}   The bilinear form
$\left( f(x), g(x) \right)_{Q,S}:=\left[ (Wf)(Q)Sg(x) \right]_{x=0}$ $f,g
\in P$ is a positive definite inner product over reals   if $n_{\psi}>0$;
$n \in N$.

\noindent {\it Proof}: $\left( s_{k}(x), s_{n}(x)
\right)_{Q,S}=[Q^{k}Ss_{n}(x)]_{x=0}=[Q^{k}q_{n}(x)]_{x=0}=n^{\underline{k}}
_{\psi}q_{n-k}(0)=n_{\psi}!\delta _{nk}$ where $\{q_{n}(x)
\}^{\infty}_{n=0}$ is the $\psi$-basic sequence of the generalised
difference-tial operator $Q$.
Let now $H=(P;(\ , \ )_{Q,S})$. We shall denote
by $(\ ,\ )_{Q,S}:H\times H \to R$
 the scalar product associated with Sheffer  $\psi$-sequence $\{s_{n}(x)
\}^{\infty}_{n=0}$. The unitary space $H =(P; (\ ,\ )_{Q,S})$ is then
completed to the unique Hilbert space $\aleph = \bar{H}$. If so then the
following Theorem is valid in $\psi$-extended and also  in $Q$-extended
and case of finite operator calculus.

\vspace{2mm}

\noindent {\bf Theorem 4.2.}       ($Q$-calculus Spectral Theorem)  [3, 16]

\noindent {\it Let $\{s_{n}(x) \}_{n=0}^{\infty}$ be generalised  Sheffer
$\psi$-sequence of $\psi$-difference-tial operator $Q$ with  $\psi$-basic
$\{ q_{n}(x)\}^{\infty}_{n=0}$ and relative to $S \in \Sigma _{Q}$. Then
there exists a unique essentially self-adjoint operator  $A_{Q,S}:\aleph
\to \aleph$
given by $A_{Q,S}=\sum_{k \geq
1}\frac{u_{k}+\hat{\nu}_{k}(x)}{(k-1)_{\psi}!}Q^{k}$ such that the spectrum
of $A_{Q,S}$ consists of $n = 0, 1, 2, 3,...$ where
$A_{Q,S}s_{n}(x)=ns_{n}(x)$. The quantities $u_{k}$ and $\hat{\nu}_{k}$ are
calculated according to:}
\[
u_{k}=-\left[(\log S)'\hat{x}_{\psi}^{-1}q_{k}(x)  \right]_{x=0}\quad
and\quad \hat{\nu}_{k}(x)=\hat{x}_{\psi}\left[ \frac{d}{dx}q_{k} \right](0).
\]

\noindent The linear map $':\Sigma _{Q} \to \Sigma _{Q}$ (as recalled by
Definition 4.4) is the Pincherle \linebreak $Q$-$\psi$-derivative. As for
the main umbral propositions [3, 16, 40] we quote here just two right after
a definition.

\vspace{2mm}

\noindent {\it Definition  4.3.}   Let $T: P\to P$ be  a  linear
operator not necessarily an element of   $\Sigma _{Q}$. If  there exist
two $\psi$-basic sequences $\{ p_{n}\}_{n \geq 0}$ and  $\{ q_{n}\}_{n \geq
0}$  of $\psi$-difference-tial operators from $\Sigma _{Q}$ such that
$Tp_{n}=q_{n},\ n \geq 0$, then we shall call  $T$ the  {\it $Q$-$\psi$-umbral
operator}.

As in [3, 40]  we may now prove
that the following theorem holds.

\vspace{2mm}

\noindent {\bf Theorem  4.3.}   {\it Let $T$  be  any $Q$-$\psi$-umbral
operator. Then for $S \in \Sigma _{Q}$ we have:

\noindent a)~the map $S \to TST^{-1}$ is an
automorphism of the algebra $\Sigma _{Q}$,

\noindent b)~if   $L$ is any   $\psi$-difference-tial operator
then  $P=TLT^{-1}$ is a $\psi$-difference-tial operator,

\noindent c)~if $S=s(Q)$ and  $L=L(Q)$  are
$\psi$-exponential formal power series in  $Q$  then $TST^{-1}=s(P)$  where
$P=TL(Q)T^{-1}$,

\noindent d)~the $Q$-$\psi$-umbral operator maps any Sheffer $Q$-$\psi$-sequence
into a Sheffer \linebreak $Q$-$\psi$-sequence.}

\noindent In order to formulate the next important  $Q$-$\psi$-umbral theorem
we need  the operator, which we shall call  the   Pincherle $Q$-$\psi$-derivative.

\vspace{2mm}

\noindent {\it Definition  4.4.}   A linear map $': \Sigma _{Q} \to
\Sigma _{Q};\quad T'=T\hat{x}_{Q}-\hat{x}_{Q}T\equiv [T,\hat{x}_{Q}] \in
\Sigma _{Q}$ is called the {\it Pincherle $Q$-$\psi$-derivative}.

As in [3, 40] we may
now prove that the following theorem holds.

\vspace{2mm}

\noindent {\bf Theorem  4.4.} {\it Let  $U$   be a $Q$-$\psi$-umbral
operator $Uq_{n}(x)=\frac{n_{\psi}!}{n!}\hat{x}_{Q}^{n}1=x^{n}$  and let
$\{q_{n}\}_{n \geq 0}$ be the $\psi$-basic sequences of
$Q$-$\psi$-difference-tial operator $L=L(Q)$. Then}
$U'=\hat{x}_{Q}U(L(Q)'-I)$.

\subsection{Miscellaneous remarks and indications of
several applications}
\setcounter{equation}{0}

\noindent {\bf A)} As announced in the Introduction we shall present now a
specific formulation of $q$-umbral calculus by Cigler [17] and
Kirschenhofer [18]. This formulation  might be related -- as noted in  [15]
-- to the so-called quantum groups [19]. Namely we shall  indicate how one
may formulate $q$-extended finite operator calculus with help of the
``quantum $q$-plane" [19]   $q$-commuting variables $A, B$:
\begin{equation} AB-qBA\equiv \left[ A, B \right]_{q}=0. \end{equation}

\noindent The idea to use ``$q$-commuting variables" goes back at least to
Cigler (1979) [17] (see formula (7), (11) in [17]) and also to
Kirchenhofer -- see [18] for further systematic development.  In [17, 18]
one defines the polynomial sequence  $\{ p_{n} \}_{0}^{\infty}$  of
$q$-binomial type by
\begin{equation}
p_{n}(A+B)\equiv \sum_{k \geq 0} \left( \begin{array}{c} n\\k \end{array}
\right)_{q} p_{k}(A)p_{n-k}(B) \quad {\rm where}\quad \left[ B, A
\right]_{q}\equiv BA-qAB=0. \end{equation}

\noindent $A$ and $B$ might be interpreted then as coordinates on quantum
$q$-plane (see [19], Chapter 4). For example $A=\hat{x}$ and $B=y\hat{Q}$
where $\hat{Q}\varphi(x)=\varphi(qx)$. If so then the following
identification takes place:
\begin{equation}
p_{n}(x+_{q}y)\equiv E^{y}(\partial _{q})p_{n}(x)=\sum _{k \geq 0}\left(
\begin{array}{c} n\\k \end{array} \right)_{q}p_{k}(x)p_{n-k}(y)=p_{n}\left(
\hat{x}+y\hat{Q}\right)1\end{equation}

\noindent Also  $q$-Sheffer  polynomials $\{s_{n}(x) \}_{n=0}^{\infty}$ are
defined equivalently (see 2.1.1 in [18])  by
\begin{equation}
s_{n}(A+B)\equiv \sum_{k \geq 0}\left( \begin{array}{c} n\\k \end{array}
\right) s_{k}(A)p_{n-k}(B) \end{equation}

\noindent where   $[B, A]_{q} \equiv BA-qAB = 0$ and  $\{
p_{n}(x)\}_{n=0}^{\infty}$ of $q$-binomial type. For example  $A=\hat{x}$
and $B=y\hat{Q}$ where $\hat{Q}\varphi(x)=\varphi(qx)$. Then the following
identification is also evident:
\begin{equation}
s_{n}(x+_{q}y)\equiv E^{y}(\partial _{q})s_{n}(x)=\sum _{k \geq 0}\left(
\begin{array}{c} n\\k \end{array} \right)_{q}
s_{k}(x)\ p_{n-k}(y)=s_{n}\left( \hat{x}+y\hat{Q}\right)1. \end{equation}

\noindent This means that one may formulate $q$-extended finite operator
calculus with help of the ``quantum $q$-plane" $q$-commuting variables $A,\
B$: $AB-qBA\equiv [A,B]_{q}=0$. A natural question then arises: is there a
$\psi$-analogue extension of quantum $q$-plane formulation? If one [15, 16]
introduces the evident natural extension of $q$-commutator as done in (see
below) then the answer given in [15, 16] is in negative. The above
identifications of polynomial sequence $\{ p_{n}\}^{\infty}_{0} $ of
$q$-binomial type and Sheffer $q$-sequences $\{s_{n}(x) \}^{\infty}_{n=0}$
fail to be extended to the more general $\psi$-case. Therefore one can not
formulate {\it that way} the $\psi$-extended finite operator calculus with
help of the ``quantum $\psi$-plane" $\hat{q}_{\psi ,Q}$-commuting variables
$A, B$: $AB - \hat{q}_{\psi, Q}BA \equiv [A, B]_{\hat{q}_{\psi},Q} = 0$. The
corresponding new objects above are defined accordingly [15].

\vspace{2mm}

\noindent {\it Definition 5.1.}   Let  $\{ p_{n} \}_{n \geq 0}$ be
the $\psi$-basic sequence of the $\psi$-difference-tial operator  $Q$. Then
the  $\hat{q}_{\psi, Q}$-operator is a liner map
\[
\hat{q}_{\psi,Q}:P\to P;\quad
\hat{q}_{\psi,Q}p_{n}=\frac{(n+1)_{\psi}-1}{n_{\psi}}p_{n};\quad n \geq 0.
\]

\vspace{2mm}
\noindent {\it Note 5.1.}   For  $Q(\partial _{\psi})=\partial _{\psi}$  the
natural notation is $\hat{q}_{\psi, id}\equiv \hat{q}_{\psi}$. For $\psi
_{n}(q)=\frac{1}{R(q^{n})!}$  and
\[ R(x)=\frac{1-x}{1-q}\
\hat{q}_{\psi,Q}\equiv \hat{q}_{R,id}\equiv \hat{q}_{R}\equiv
\hat{q}_{q,id}\equiv \hat{q}_{q}\equiv \hat{q}\]

\noindent and
\[\hat{q}_{\psi,Q}x^{n}=q^{n}x^{n} . \]

\vspace{2mm}

\noindent {\it Definition 5.2.}   Let  $A$ and $B$ be linear operators
acting on $P$;  $A:P\to P;$ $B:P\to P$. Then
$AB-\hat{q}_{\psi,Q}BA\equiv[A,B]_{\hat{q}_{\psi,Q}}$ is called
$\hat{q}_{\psi,Q}$-mutator of $A$ and $B$ operators.

Naturally  $Q\hat{x}_{Q}-\hat{q}_{\psi ,Q}\hat{x}_{Q}Q\equiv [Q,
\hat{x}_{Q}]_{\hat{q}_{\psi,Q}}={\rm id}$. This is easily verified in the
$\psi$-basic sequence of the operator $Q$. In conclusion the case of
$q$-extended finite operator calculus -- or $q$-calculus in short -- is
fairly enough distinguished by the Cigler-Kirchenhofer approach and the
quantum plane notion among the infinite variety of  $Q$-$\psi$-umbral
calculi.

\vspace{3mm}

\noindent {\bf B)}
We end our exposition of the fundamentals of   $\psi$-umbral calculus with
indications of several {\it examples of identities} resulting
straightforwardly from the fact that  dual pairs $\{ Q, \hat{x}_{Q} \}$  --
(and in particular $\{Q(\partial _{\psi}), \hat{x}_{Q(\partial _{\psi})}
\}$ including any dual pair $\{ \partial _{\psi}, \hat{x}_{\psi} \}$
determined  by any choice of   admissible  sequence  $\psi$) --   provide
us with distinct $Q$-$\psi$-representations of  GHW algebra as seen from
$End (P) = [\{\partial _{\psi},\hat{x}_{\psi} \}]=[\{Q(\partial
_{\psi}),\hat{x}_{Q(\partial _{\psi})}\}]=[\{ Q, \hat{x}_{Q} \}]$.
\begin{enumerate}
\item Let $\{a_{n}(x)\}_{n \geq 0}$ be any Appel  $\psi$-sequence. Then
(compare (4) in [41]) we have
\begin{equation}
\hat{x}_{Q}\sum_{m=0}^{n} a_{m}(Q)\frac{\hat{x}_{Q}^{m}}{m_{\psi}!}\equiv
a_{n}(Q)\frac{\hat{x}_{Q}^{n+1}}{m_{\psi}!}.
\end{equation}

Also other identities from [41--42]  apply specifically to dual pairs
$\{ Q,\hat{x}_{Q} \}$. For example:

\item Compare  (6)  and  (7)  in [42]
\begin{equation}
\left( Q\hat{x}_{Q} Q \right)^{n}=Q^{n}\hat{x}_{Q}^{n}Q^{n};\quad \left(
\hat{x}_{Q}Q\hat{x}_{Q}\right)^{n}=\hat{x}_{Q}^{n}Q^{n}\hat{x}_{Q}^{n}.
\end{equation}

\item Similarly identities from [43]  apply to dual pairs $\{Q, \hat{x}_{Q}
\}$. For example one has for $\hat{A}^{\underline{n_{\psi}}}=\hat{A}\left(
\hat{A}-1_{\psi}{\rm id} \right)\left( \hat{A}-2_{\psi}{\rm id} \right)...\left(
\hat{A}-(n-1)_{\psi}{\rm id} \right)$  (compare with  (8)  in [43]) the identity
\begin{equation}
\hat{x}_{Q}^{n}Q^{n}\left[ f(\hat{x}_{Q}) \right]=\left(\hat{x}_{Q}Q
\right)^{\underline{n_{\psi}}}\left[ f(\hat{x}_{Q}) \right]. \end{equation}

\item Lagrange inversion formula and more general formula from the theorem due to Viskov
(in [29] point 18$^{\circ}$) also apply automatically to dual pairs $\{
Q,\hat{x}_{Q}\}$ under evident replacements.
\end{enumerate}






\vspace{0.5cm}
\noindent {\erm  Institute of Computer Science}

\noindent{\erm Bia\l ystok University }

\noindent{\erm
Sosnowa 64, PL-15-887 Bia\l ystok}

\noindent {\erm Poland}

\vspace{0.5cm}
\noindent Presented by Julian \L awrynowicz at the Session of the
Mathematical-Physical Commission of the \L \'od\'z Society of Sciences and
Arts on September 26, 2002

\vspace{0.5cm}
\noindent {\bf PROSTE CHARAKTERYZACJE $\psi$-WIELOMIAN\'OW SHEFFERA
\linebreak I WYNIKAJ\c{A}CE ST\c{A}D IMPLIKACJE DOTYCZ\c{A}CE RACHUNKU
CI\c{A}G\'OW} \vspace{0.2cm}

\noindent {\small S t r e s z c z e n i e}

{\small ,,Rachunek ci\c{a}g\'ow" zapocz\c{a}tkowa\l a publikacja Warda
(1936) sugeruj\c{a}ca mo\.zliwy zakres rozszerze\'n rachunku operatorowego
Roty-Mullina, rozwa\.zanych po Wardzie przez wielu \linebreak autor\'ow.
Dla wygody opracowane p\'o\'zniej warianty rachunku ci\c{a}g\'ow Warda
b\c{e}dziemy nazywali $\psi$-rachunkiem. Oznaczenia u\.zywane przez Warda,
Viskova, Markowsky'ego i Romana s\c{a} w niniejszej pracy zastosowane w
uzgodnieniu z oznaczeniami Roty. W ten spos\'ob \linebreak $\psi$-rachunek staje
si\c{e} cz\c{e}\'sciowo niemal automatycznym rozszerzeniem sko\'nczonego
rachunku operator\'ow. $\psi$-rozszerzenie opiera si\c{e} na poj\c{e}ciu
niezmienniczo\'sci operator\'ow wzgl\c{e}dem
$\partial_{\psi}$-przesuni\c{e}cia. Jednocze\'snie omawiany rachunek jest
przyk\l adem algebraizacji analizy -- tutaj algebraizacji ograniczonej do
algebry szereg\'ow formalnych. Dogodno\'s\'c u\.zytych oznacze\'n wyra\.za
si\c{e} mi\c{e}dzy innymi przez mo\.zliwo\'s\'c \l atwego udowodnienia
pewnych twierdze\'n o charakteryzacji $\psi$-wielomian\'ow Sheffera jak
r\'ownie\.z twierdzenia spektralnego. Wyniki $\psi$-rachunku mog\c{a} by\'c
rozszerzone na ,,rachunek {\eit Q}-umbralny Markowsky'ego", gdzie {\eit Q}
oznacza uog\'olniony operator r\'o\.znicy (nie musi on by\'c niezmienniczy
wzgl\c{e}dem $\partial_{\psi}$-przesuni\c{e}- \linebreak cia), tj. operator
obni\.zaj\c{a}cy  stopie\'n dowolnego wielomianu o jeden.}



\end{document}